\documentclass[authoryear,3p,preprint,11pt]{elsarticle}

\usepackage[utf8]{inputenc}
\usepackage{amsmath}
\usepackage[english]{babel}
\usepackage{amssymb}
\usepackage{amsfonts}
\usepackage{bm}
\usepackage{graphicx}
\usepackage{amsmath}
\usepackage{amsthm}
\usepackage{xcolor}
\usepackage{parskip}
\usepackage{multicol}
\usepackage{float}
\usepackage{rotating}
\usepackage{subcaption}
\usepackage{rotating}

\usepackage{chngcntr}
\usepackage{apptools}
\AtAppendix{\counterwithin{lem}{subsection}}
\AtAppendix{\counterwithin{proposition}{subsection}}
\AtAppendix{\counterwithin{equation}{subsection}}

 \theoremstyle{definition}

\newtheorem{proposition}{Proposition}[section]
\newtheorem{mydef}{Definition}[section]

\newtheorem{remark}{Remark}[section]

\usepackage[colorlinks=true,linkcolor=red, citecolor=blue, urlcolor=red]{hyperref}

\usepackage{amssymb}

 \usepackage{lineno}

\journal{}

\usepackage[symbol]{footmisc}

\begin{document}


\begin{frontmatter}


\title{Estimating covariance functions of multivariate skew-Gaussian random fields on the sphere}


\author[label1]{Alegr\'ia A. \footnote{Corresponding author. Email: alfredo.alegria.jimenez@gmail.com}}
\author[label2]{Caro S.}
\author[label3]{Bevilacqua M.}
\author[label1]{Porcu E.}
\author[label1]{Clarke J.}
\address[label1]{Departamento de Matem\'atica, Universidad T\'ecnica Federico Santa Mar\'ia, Valpara\'iso, Chile}
\address[label2]{Departamento de Matem\'atica y Computaci\'on, Universidad de Santiago, Santiago, Chile}
\address[label3]{Instituto de Estad\'istica, Universidad de Valpara\'iso, Valpara\'iso, Chile}

\begin{abstract}
 This paper considers a multivariate  spatial random field, with each component having univariate marginal distributions of the skew-Gaussian type.  We assume that the field is defined spatially on the unit sphere embedded in $\mathbb{R}^3$, allowing for  modeling  data available over large portions of planet Earth.   This model admits explicit expressions for the  marginal and cross covariances. However,  the $n$-dimensional distributions of the field are difficult to evaluate, because it requires the sum of $2^n$ terms involving the cumulative and probability density functions of a $n$-dimensional Gaussian distribution. Since  in this case  inference based on the full likelihood is computationally unfeasible, we propose a composite likelihood approach based on pairs of spatial observations.   This last being possible thanks to the fact that we have a closed form expression for the bivariate distribution. We illustrate the effectiveness of the method through simulation experiments and the analysis of a real data set of minimum and maximum surface air temperatures. 
\end{abstract}

\begin{keyword}
Composite likelihood \sep Geodesic distance \sep Global data

\end{keyword}

\end{frontmatter}



\section{Introduction}

The Gaussian assumption is an appealing option  to model  spatial data. First, the Gaussian distribution is completely characterized by its first two moments. Another interesting property is the tractability of the Gaussian distribution under linear combinations and conditioning. However, in many geostatistical applications, including oceanography, environment and the study of natural resources,  the Gaussian framework    is unrealistic, because the observed data have  different features, such as, positivity, skewness or heavy tails, among others.

Transformations of Gaussian random fields (RFs)  is  the most common alternative to model non-Gaussian fields. Consider a spatial domain $\mathcal{D}$ and   $\{Z(\bm{s}), \bm{s}\in \mathcal{D}\}$ defined as
  $Z(\bm{s})=\varphi(Y(\bm{s}))$,  
where  $\varphi$ is a real-valued mapping and $\{Y(\bm{s}), \bm{s}\in \mathcal{D}\}$ is Gaussian.
Apparently the finite-dimensional distributions of $Z$ depend on the choice of  $\varphi$.  In some cases, such mapping is  one-to-one  and admits an inverse simplifying the analysis. In this class, we can highlight the log-Gaussian RFs, which are generated as a particular example of the Box-Cox transformation  (see \citealp{deoliveira:1997}).  However, a one-to-one transformation is  not appropriate in general. For instance, \cite{stein:92} considers a truncated Gaussian RF, taking $\varphi$ as an indicator function, in order   to model data sets with a given percentage of  zeros (for instance, precipitations). On the other hand, wrapped-Gaussian RFs have been introduced in the literature for  modeling directional spatial data, arising in the study of wave and wind directions (see \citealp{jona-lasinio2012}).   In addition, \cite{xu2016tukey} propose a flexible class of fields named the Tukey g-and-h RFs.

Another approach consists in taking advantage of the stochastic representations of random variables. For instance,  \cite{ma:2009}  considers a general approach to construct elliptically contoured RFs through mixtures of Gaussian fields. These models have an explicit covariance structure and  allow a wide range of finite-dimensional distributions.  Other related works have been developed by \cite{du2012hyperbolic}, \cite{ma2013k} and \cite{ma2013student} including Hyperbolic, K and Student's t distributed fields.
 Moreover,  \cite{kim:mallick:2004}, \cite{gualtierotti:2005} and \cite{allard:naveau:2007}  have  introduced skew-Gaussian RFs for  modeling  data with skewed distributions. 
   However, \cite{Minozzo:2012} and \cite{Genton:2012} show that all these models  are not valid because they cannot be identified with probability one using a single realization, i.e.,  in practice it is impossible to make inference on the basis of these models. Such results do not prevent the existence of  RFs having univariate marginal distributions belonging to a given  family. 
   
   In this paper, we consider  multivariate stationary  RFs, where  each component has a univariate marginal distribution of the skew-Gaussian type. We follow the representation proposed in the univariate case  by \cite{zhang:el-shaarawi:2010} and extend it to the multivariate case. This construction allows for modeling   data with different degrees of  skewness   as well as explicit expressions for the covariance function.   Estimation methods for this model are still unexplored. Maximum likelihood is certainly  a useful tool, but it is  impracticable, because the full likelihood  does not have a simple form.  Indeed,  if $n$ is the number of observations, it can be explicitly expressed as the   sum of $2^{n}$ terms depending on the probability density function (pdf) and the cumulative distribution function (cdf) of the $n$-variate Gaussian distribution. Direct maximization of the likelihood seems intractable from a computational and analytical point of view.    
   \cite{zhang:el-shaarawi:2010} consider the EM algorithm to perform inference on the skew-Gaussian model.  However, the iterations of the EM algorithm are difficult to evaluate because each step requires Monte Carlo simulations of a non-trivial conditional  expectation.     
 On the other hand, composite likelihood (CL) is an estimation procedure (\citealp{Lindsay:1988}; \citealp{Varin:Reid:Firth:2011}; \citealp{Cox:Reid:2004}) based on the likelihood of marginal or conditional events. CL methods are an attractive option when the full likelihood is difficult to write  and/or when the data sets are large. This approach has been used in several spatial and space-time contexts, mainly in the Gaussian case (\citealp{Vecchia:1988}; \citealp{Curriero:Lele:1999};  \citealp{Stein:Chi:Welty:2004};  \citealp{Bevilacqua:Gaetan:Mateu:Porcu:2012};  \citealp{Bevilacqua:Gaetan:2014}; \citealp{Bevilacqua2016}). Outside the Gaussian scenario, \cite{Heagerty:Lele:1998} propose  CL inference for  binary spatial data.  Moreover, \cite{Padoan:Ribatet:Sisson:2010} and \cite{sang:2014} have used  CL methods for the estimation of max-stable fields, whereas \cite{doi:10.1080/00949655.2016.1162309} consider a  truncated CL approach for wrapped-Gaussian fields. 
 
 The implementation of the  CL method on multivariate skew-Gaussian fields is still unexplored. Our goal consists in developing a CL  approach  based on pairs of observations for a multivariate skew-Gaussian RF.    Our contribution provides a fast and accurate tool to make inference on skewed data.  The main ingredient of the pairwise CL method is the characterization of the bivariate distributions of the RF, that is, we derive a closed form expression for the joint distribution between two correlated  skew-Gaussian random variables (possibly with different means, variances and degrees of skewness).

 In addition, in order to work with data collected over the whole planet Earth, we consider the spatial domain as the unit sphere $\mathcal{D}=\mathbb{S}^2:=\{\bm{s}\in\mathbb{R}^3: \|\bm{s}\|=1\}$, where $\|\cdot\|$ denotes the Euclidean distance. We refer the reader to \cite{marinucci2011random} for a more detailed study about RFs on spheres.  An important implication is that the covariance function depends on a different metric, called geodesic distance. In general, covariance models valid on Euclidean spaces are not valid on the sphere, and we refer the reader to \cite{gneiting2013} for a overview of the problem. 
 
The paper is organized as follows. In Section 2, the multivariate skew-Gaussian RF is introduced. The bivariate distributions of the skew-Gaussian field and the  CL approach are discussed in Section 3. In section 4,  simulation experiments are developed.  Section 5 contains a real application for a bivariate data set of minimum and maximum surface air temperatures.  Finally, Section 6 provides a brief discussion.

\section{Skew-Gaussian RFs on the unit sphere}

In this section, we introduce a skew-Gaussian model  generated as a mixture of two latent  Gaussian RFs.  Such construction is based on the stochastic representation of skew-Gaussian random variables.  Let $X$ and $Y$ be two independent standard Gaussian random variables and $-1\leq \delta \leq 1$. Then, the distribution of  
\begin{equation}\label{stoch_repr}
Z=\delta |X| + \sqrt{1-\delta^2}Y
\end{equation}
 is called skew-Gaussian, with pdf $2\phi(z)\Phi(\alpha z)$, for $z\in\mathbb{R}$, where $\alpha=\delta/\sqrt{1-\delta^2}$. Here, $\phi(\cdot)$ and $\Phi(\cdot)$ denote the pdf and cdf of the standard Gaussian distribution.  If $\delta >0$, we say that $Z$ is right-skewed, whereas for $\delta<0$, $Z$ is left-skewed. Of course, for $\delta=0$ we have the Gaussian case.
For a detailed study of the skew-Gaussian distribution, we refer the reader to \cite{azzalini1985class,azzalini1986further}, \cite{azzalini1996multivariate}, \cite{azzalini1999statistical}, \cite{arellano:2006} and \cite{azzalini2013skew}. 

We work with the spatial counterpart of  Equation (\ref{stoch_repr}). Let   $\{  \bm{X}(\bm{s}) = (X_1(\bm{s}),\hdots,X_m(\bm{s}))^\top: \bm{s}\in\mathbb{S}^2\}$ and $\{ \bm{Y}(\bm{s}) = (Y_1(\bm{s}),\hdots,Y_m(\bm{s}))^\top: \bm{s}\in\mathbb{S}^2\}$ be two   stationary multivariate  Gaussian  RFs  defined on $\mathbb{S}^2$.  Here, $m\in\mathbb{N}$ denotes the number of components of the fields and $\top$ is the transpose operator. In addition, we  assume that $\bm{X}(\bm{s})$ and $\bm{Y}(\bm{s})$ are independent, with components having zero mean and unit variance.
 In the spherical framework, the covariances are given  in terms of  the geodesic  distance, defined by $$\theta :=\theta(\bm{s}_1,\bm{s}_2) = \arccos (\bm{s}_1^\top \bm{s}_2) \in[0,\pi], \qquad \bm{s}_1,\bm{s}_2\in\mathbb{S}^2,$$ which is the most natural metric on the spherical surface.
 Therefore, we suppose that there exists two matrix-valued  mappings $r^x,r^y: [0,\pi]\rightarrow \mathbb{R}^{m\times m}$ such that 
\begin{equation*}
\text{cov}\{X_i(\bm{s}),X_j(\bm{s}')\}   = r^x_{ij}(\theta(\bm{s},\bm{s}'))  \text{  and  }   \text{cov}\{Y_i(\bm{s}),Y_j(\bm{s}')\}  = r^y_{ij}(\theta(\bm{s},\bm{s}')),
\end{equation*}
for all $\bm{s},\bm{s}'\in\mathbb{S}^2$  and $i,j=1,\hdots,m$. In such case, we say that the covariance  function is spherically isotropic.  In the univariate case ($m=1$), \cite{gneiting2013} establishes that some classical covariances such as the Cauchy, Mat\'ern, Askey and Spherical models, given in the classical literature in terms of Euclidean metrics, can be coupled with the geodesic distance under specific constraints for the parameters.  Furthermore, \cite{PBG16} propose several  covariance models for space-time and multivariate  RFs on spherical spatial domains.

Next, we define a multivariate spatial RF with each component being skew-Gaussian distributed, according to Equation (\ref{stoch_repr}). It is a multivariate extension of the univariate  skew-Gaussian field proposed by \cite{zhang:el-shaarawi:2010}. This model allows  different means, variances and  skewness in the components of the RF.
\begin{mydef}
A  multivariate field, $\{  \bm{Z}(\bm{s}) =(Z_1(\bm{s}),\hdots,Z_m(\bm{s}))^\top: \bm{s}\in\mathbb{S}^2\}$, with components having skew-Gaussian marginal distributions, can be defined through
\begin{equation}
\label{model_skew}
Z_i(\bm{s}) = \mu_i +  \eta_i |X_i(\bm{s})| + \sigma_i Y_i(\bm{s}) , \qquad \bm{s}\in\mathbb{S}^2, \qquad i=1,\hdots,m,
\end{equation}
where $\mu_i,\eta_i \in \mathbb{R}$ and  $\sigma_i \in \mathbb{R}_+$. Note that  $(Z_i(\bm{s})-\mu_i)/\sqrt{\eta_i^2+\sigma_i^2}$, for all $\bm{s}\in\mathbb{S}^2$, follows a skew-Gaussian distribution with pdf given by
$ f_{Z_i}(z) = 2 \phi\left( z  \right) \Phi\left( (\eta_i/\sigma_i) z \right).$
\end{mydef}

\begin{remark}
Recent literature considers the latent fields $X_i(\bm{s})$, $i=1,\hdots,m$,  as single random variables $X_i$,   being constants along  the spatial domain. However, this approach has apparent identifiability problems,  since in practice we only  work with one realization and there is no information about the variability of $X_i$. Thus, this approach only produces  a shift effect in the model. These considerations are studied by \cite{Minozzo:2012} and \cite{Genton:2012}. 
\end{remark}

Direct application of the results given by \cite{zhang:el-shaarawi:2010} provides the following proposition.
 \begin{proposition}
The field $\bm{Z}(\bm{s})$ defined through Equation (\ref{model_skew})  is   stationary  with expectations
\begin{equation*}
\mathbb{E}(Z_i(\bm{s}) ) =  \mu_i + \eta_i \sqrt{\frac{2}{\pi}}, \qquad \bm{s}\in\mathbb{S}^2,
\end{equation*}
and covariances
\begin{equation}
\label{cov_skew}
C_{ij}(\theta) := \text{cov}\{ Z_i(\bm{s}), Z_j(\bm{s}') \} =    \frac{2\eta_i\eta_j}{\pi}  g(r_{ij}^x(\theta)) + \sigma_i\sigma_j  r_{ij}^y(\theta),
\end{equation}
for all $i,j=1,\hdots,m$, where   $\theta= \theta(\bm{s},\bm{s}')$ and $g(t)=\sqrt{1-t^2}+t\arcsin(t) -1$, for $|t|\leq 1$.
\end{proposition}
The proof is omitted because it is obtained by using the same arguments as in \cite{zhang:el-shaarawi:2010}.

\section{Composite likelihood estimation}

\subsection{General framework}

We first introduce the CL approach from a general point of view. CL methods \citep{Lindsay:1988} are likelihood  approximations for dealing with large data sets. In addition, in the last years,  these techniques have been used to study data with  intractable analytical expressions for the full likelihood. The objective function for CL methods is constructed through  the likelihood of marginal or conditional events. Formally, let $f(\bm{z};\bm{\lambda})$ be the pdf of a  $n$-dimensional random vector, where $\bm{\lambda}\in \Lambda \subset \mathbb{R}^{p}$ is an unknown parameter vector, and $\Lambda$ is the parametric space. We denote by $\{ \mathcal{A}_{1} ,..., \mathcal{A}_{K} \}$ a set of marginal or conditional events with associated likelihoods $\mathcal{L}_{k}(\bm{\lambda};\bm{z})$. Then, the objective function for the composite likelihood method is defined as the weighted product
  \begin{equation*}
  \mathcal{L}_{C}(\bm{\lambda};\bm{z}) := \displaystyle \prod_{k=1}^{K} \mathcal{L}_{k}(\bm{\lambda};\bm{z})^{w_{k}},
 \end{equation*}
where the non-negative weights $w_{k}$ must be chosen according to an appropriate criterion. In principle, the weights can  improve  the statistical and/or computational efficiency of the estimation. We use the notation $\ell_k(\bm{\lambda};\bm{z})=\log \mathcal{L}_k(\bm{\lambda};\bm{z})$, thus the log composite likelihood is given by
\begin{equation*}
 \text{CL}(\bm{\lambda};\bm{z}) := \displaystyle \sum_{k=1}^{K} \ell_{k}(\bm{\lambda};\bm{z}) w_{k}.
\end{equation*}
The maximum CL estimator is defined as $$\widehat{\bm{\lambda}}_{n}=\text{argmax}_{\bm{\lambda}\in\Lambda}\text{CL}(\bm{\lambda};\bm{z}).$$
   By construction, the composite score $\nabla \text{CL}(\bm{\lambda})$ is an unbiased estimating equation, i.e., $\mathbb{E}(\nabla \text{CL}(\bm{\lambda}))=\bm{0} \in \mathbb{R}^p$. This is an  appealing property of CL methods, since   it is a   first order likelihood property.  On the other hand, the second order properties are related to the  Godambe information matrix, defined as $$G_{n}(\bm{\lambda})=H_{n}(\bm{\lambda})J_{n}(\bm{\lambda})^{-1}H_{n}(\bm{\lambda})^\top,$$
where $H_{n}(\bm{\lambda})=-\mathbb{E}(\nabla^{2}\text{CL}(\bm{\lambda};\bm{z}))$ and $J_{n}(\bm{\lambda})=\mathbb{E}(\nabla \text{CL}(\bm{\lambda};\bm{z}) \nabla \text{CL}(\bm{\lambda};\bm{z})^\top)$. The inverse of $G_{n}(\bm{\lambda})$ is an approximation of the asymptotic variance of the CL estimator. Under increasing domain and regularity assumptions, CL estimates are consistent and asymptotically Gaussian.

\subsection{Pairwise CL approach for the  multivariate skew-Gaussian model}

We now develop a CL method based on pairs of observations for the multivariate skew-Gaussian RF.  We consider the $m$-variate field, $\{ \bm{Z}(\bm{s})=(Z_1(\bm{s}),\hdots,Z_m(\bm{s}))^\top, \bm{s}\in \mathbb{S}^2 \}$, defined in  (\ref{model_skew}), and a   realization of $\bm{Z}(\bm{s})$ at $n$ spatial locations, namely,  $(\bm{Z}(\bm{s}_{1})^\top,\hdots,\bm{Z}(\bm{s}_{n})^\top)^\top$. Then, we define all possible   pairs   $\bm{Z}_{ijkl}=(Z_i(\bm{s}_{k}),Z_j(\bm{s}_{l}))^\top$  with associated log likelihood $\ell _{ijkl}(\bm{\lambda})$, where $\bm{\lambda}\in \mathbb{R}^{p}$ is the parameter vector. Therefore, the corresponding  log composite likelihood equation is defined by (see \citealp{Bevilacqua2016})
\begin{equation*}
\label{cl}
\text{CL}(\bm{\lambda}) = \sum_{i=1}^{m}\sum_{k=1}^{n-1} \sum_{l=k+1}^{n}  \omega_{iikl} \ell_{iikl}(\bm{\lambda
}) + \sum_{i=1}^{m-1}\sum_{j=i+1}^{m}\sum_{k=1}^{n} \sum_{l=1}^{n}  \omega_{ijkl} \ell_{ijkl}(\bm{\lambda}).
\end{equation*}
Note that the order of computation of the method is $\mathcal{O}\{ mn(n-1)/2 + m(m-1)n^2/2  \}$. Throughout, we consider a cut-off weight function (0/1 weights):  $w_{ijkl}= 1$ if  $\theta(\bm{s}_{k} , \bm{s}_{l}) \leq d_{ij}$, and 0 otherwise, for some cut-off distances $d_{ij}$. This choice has apparent computational advantages. Moreover, it can improve the efficiency as it has been shown in \cite{Joe:Lee:2009}, \cite{Davis:Yau:2011} and \cite{Bevilacqua:Gaetan:Mateu:Porcu:2012}. The intuition behind this approach is that the correlations between pairs of distant observations are often nearly zero. Therefore, using all possible pairs can generate a  loss of efficiency, since redundant pairs  can produce bias in the results.

  From now on, we use the notation 
\begin{equation}\label{corr_matrix}
\Omega(r)=\begin{pmatrix}
1 & r\\
r & 1
\end{pmatrix}.
\end{equation}
The following result characterizes    the pairwise distributions for the multivariate skew-Gaussian RF. We suppose that $\eta_i\neq 0$, for all $i=1,\hdots,m$. The case with zero skewness is reduced to the Gaussian scenario. 
\begin{proposition}
\label{biv_log_lik}
Consider  two sites $\bm{s}_k,\bm{s}_l\in\mathbb{S}^2$ and  $\theta=\theta(\bm{s}_k,\bm{s}_l)$. The log likelihood associated to the pair $\bm{Z}_{ijkl}=(Z_i(\bm{s}_{k}),Z_j(\bm{s}_{l}))^\top$ in the multivariate skew-Gaussian model (\ref{model_skew}) is given by
\begin{equation}
  \label{lik_par}
\ell _{ijkl}(\bm{\lambda})  = \log\left( \displaystyle  2  \sum_{t=1}^{2}    \phi_{2} (\bm{Z}_{ijkl}-\bm{\mu};A_{t} ) \Phi_{2} ( L_{t};    B_{t}  )  \right)  
\end{equation}
where  $\phi_{2} (\bm{y};\Sigma )$ denotes the  bivariate Gaussian density function with zero mean and covariance matrix $\Sigma$. Similarly,  $\Phi_{2} ( \bm{l}; \Sigma)$ denotes the corresponding Gaussian cdf.  Here, 
\begin{eqnarray*}
\bm{\mu}  &    =  & (\mu_i,\mu_j)^\top,\\
A_{t}  &   =  &  \Omega_2 + \Upsilon ^{-1} \Omega[(-1)^{t} r_{ij}^x(\theta)] \Upsilon ^{-1},\\
B_{t}  &   =  &   ( [\Upsilon \Omega_2 \Upsilon]^{-1} + \Omega[(-1)^{t} r_{ij}^x(\theta)]^{-1} )^{-1}  ,\\
L_{t} &   =  & \bigg[  I_2 + \Upsilon \Omega_2 \Upsilon \Omega[(-1)^{t} r_{ij}^x(\theta)]^{-1}   \bigg]^{-1} \Upsilon (\bm{Z}_{ijkl}-\bm{\mu}), 
\end{eqnarray*}
 where  $\Upsilon= \text{diag}\{1/\eta_i, 1/\eta_j\}$,  $I_2$ is the identity matrix of order $(2\times 2)$ and $$\Omega_2 = \begin{pmatrix}  \sigma_i^2  &  \sigma_i\sigma_j  \\ \sigma_i\sigma_j &  \sigma_j^2 \end{pmatrix} \circ \Omega(r_{ij}^y(\theta)),$$
where $\circ$ denotes the Hadamard product.

We have deduced a closed form expression for the bivariate distributions of the  field. Note that the correlation function $r_{ij}^{x}$ alternates its sign in each element of the sum. Evaluation of Equation (\ref{lik_par}) requires the numerical calculation   of the bivariate Gaussian cdf.  The proof of Proposition \ref{biv_log_lik} is deferred to  Appendix A.
 
\end{proposition}

\section{Simulation study}

This section assesses  through simulation experiments the statistical and computational performance of the pairwise CL method. We pay attention to bivariate   ($m=2$) skew-Gaussian RFs on $\mathbb{S}^2$.

\subsection{Parameterization}

We believe that there are no strong arguments to consider different correlation structures $r^x(\cdot)$ and $r^y(\cdot)$ for the latent RFs. For example, the smoothness of the skew-Gaussian RF  is the same as the smoothness of  the roughest latent Gaussian field. Moreover, if both latent correlations are compactly supported, thus the covariance generated has also compact support.   We thus consider  latent fields belonging to the same parametric family of correlation functions, 
$$r_{ij}^x(\theta) =  \rho_{ij}^x r(\theta;c_{ij}^x),   \qquad    r_{ij}^y(\theta) = \rho_{ij}^y r(\theta;c_{ij}^y),   \qquad i,j=1,2, \qquad \theta\in[0,\pi],$$ 
where  $\rho_{ii}^x=1$,  $\rho_{12}^x = \rho_{21}^x$,  $|\rho_{12}^x|\leq 1$ and $c_{ij}^x>0$ (and similar conditions for $\rho_{ij}^y$ and $c_{ij}^y$), with mapping $\theta \mapsto r(\theta;c)$ being any univariate correlation function on the sphere (see \citealp{gneiting2013}).  

The  particular choice of $r(\cdot;\cdot)$  produces additional restrictions on the parameters. Throughout,  we  work with 
\begin{equation}\label{exp_corr}
r(\theta;c) = \exp\left(  -\frac{3\theta}{c}\right), \qquad \theta\in[0,\pi],
\end{equation}
or
\begin{equation}\label{askey_corr}
r(\theta;c) = \left(1-  \frac{\theta}{c} \right)_{+}^4, \qquad \theta\in[0,\pi],
\end{equation}
where $c>0$ is  a scale parameter and $(a)_+ = \max\{0,a\}$. Mappings (\ref{exp_corr}) and (\ref{askey_corr}) are known as Exponential and Askey models, respectively.  The former decreases   exponentially to zero and it takes values less than $0.05$ for $\theta > c$, whereas  the second is compactly supported, that is, it is identically equal to zero beyond the cut-off distance $c$.  Explicit parametric conditions  for the  validity of the bivariate Exponential model  are provided by  \cite{PBG16}. On the other hand,   Appendix B   illustrates a construction principle  that justify the use of  bivariate Askey models on  spheres.  

\begin{figure}
\centering{
\includegraphics[scale=0.35]{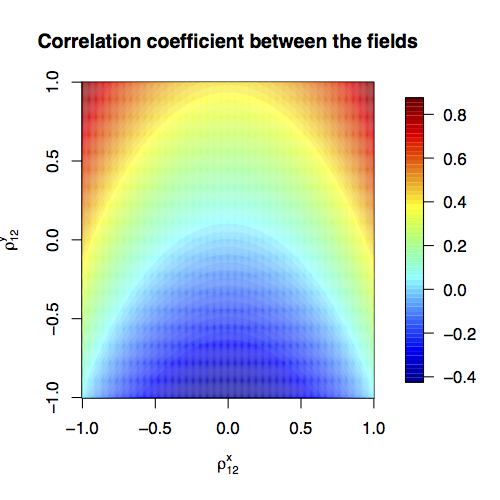} \includegraphics[scale=0.35]{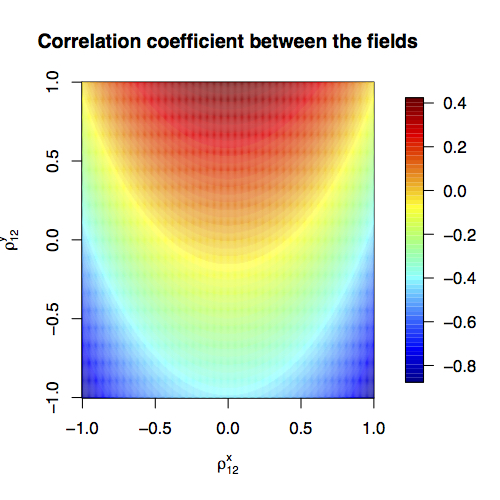}
} 
\caption{Collocated correlation coefficient  between the components of a bivariate skew-Gaussian RF in terms of $\rho_{ij}^x$ and $\rho_{ij}^y$. We consider  $\sigma_1^2 = \sigma_2^2 =1$    and  two scenarios for the skewness parameters: $(\eta_1,\eta_2)=(1,3)$ (left) and $(\eta_1,\eta_2)=(1,-3)$ (right).}
\label{contours0}
\end{figure}

An interesting property is that the   collocated correlation coefficient between the components of a bivariate skew-Gaussian RF,  $C_{12}(0) / \sqrt{C_{11}(0)C_{22}(0)}$, with $C_{12}(\cdot)$ defined in (\ref{cov_skew}),   depends on the majority of the model parameters. Figure \ref{contours0} shows the behavior of this coefficient in terms of the latent correlation coefficients $\rho_{ij}^x$ and $\rho_{ij}^y$, with $\sigma_1^2 = \sigma_2^2 =1$    and  under two different settings for the skewness parameters: $(\eta_1,\eta_2)=(1,3)$ and $(\eta_1,\eta_2)=(1,-3)$.  Note that the correlation between left-skewed and right-skewed fields has a more restrictive upper bound and this case admits strong negative correlations.   The following studies  are based on the parsimonious parameterizations $\rho_{12} := \rho_{12}^x = \rho_{12}^y$,  $c_{ij} := c_{ij}^x = c_{ij}^y$ and $c_{12}=(c_{11}+c_{22})/2$.   The parameter vector is given by $\bm{\lambda} = (\sigma_1^2, \sigma_2^2, \eta_1, \eta_2, c_{11}, c_{22}, \rho_{12}, \mu_1, \mu_2)^\top$.   In addition,  this parameterization  avoids identifiability problems.

\begin{figure}
\centering{
\includegraphics[scale=0.3]{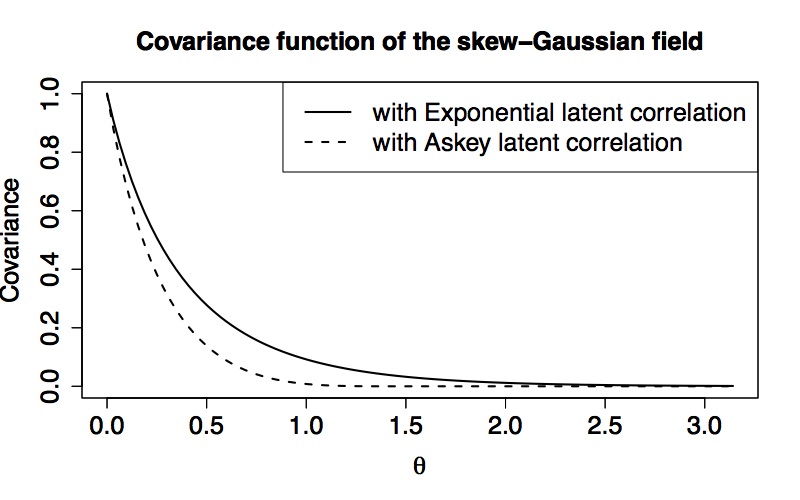}  
}
\caption{Covariance function associated to the skew-Gaussian RF, with latent correlations of Exponential (solid line) and Askey  (dashed line) types. }
\label{curves}
\end{figure}

\begin{figure}
\centering{
\includegraphics[scale=0.22]{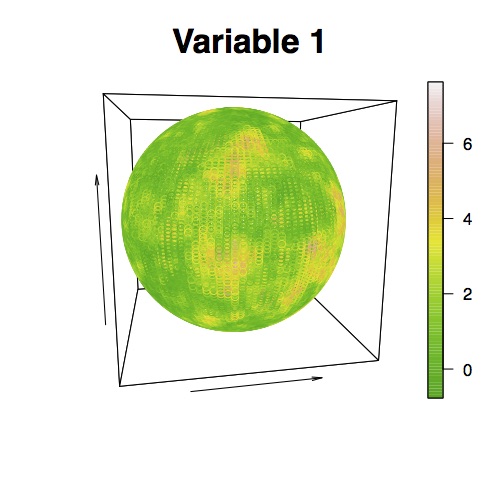}  \includegraphics[scale=0.22]{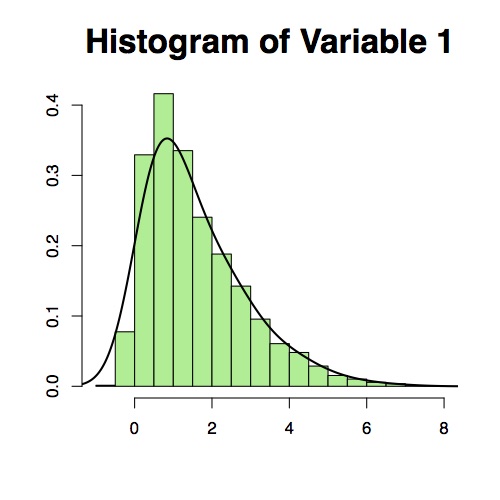}    \includegraphics[scale=0.22]{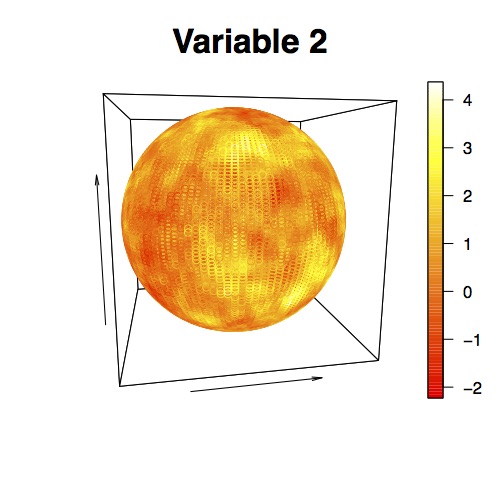}  \includegraphics[scale=0.22]{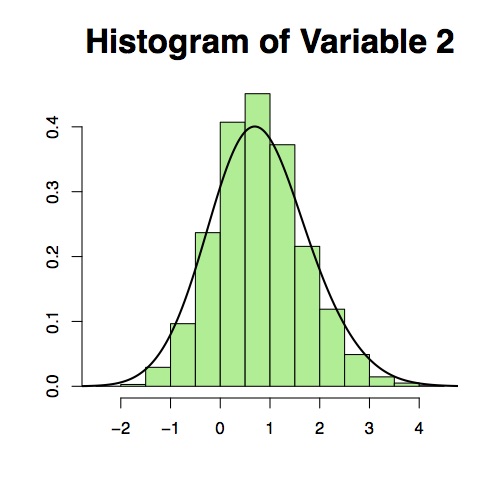} 
}
\caption{Bivariate simulation from the skew-Gaussian model with latent fields having  Exponential correlation functions. Both variables are right-skewed and the empirical  collocated correlation coefficient  is approximately 0.7.}
\label{simulation_exp}
\end{figure}

Figure \ref{curves}  shows the  covariance $C_{12}(\theta)$, in  Equation (\ref{cov_skew}), generated from         the latent correlation functions (\ref{exp_corr}) and (\ref{askey_corr}).   The skew-Gaussian RF preserves the correlation shape of the latent fields.     Figure \ref{simulation_exp} depicts a  bivariate realization of a skew-Gaussian RF, over 15000 spatial locations, with latent fields having  Exponential  correlations. We have simulated using Cholesky decomposition with $\sigma_1^2 = 0.1$, $\sigma_2^2=0.5$, $\eta_1 = 2$, $\eta_2=1$, $\mu_1=\mu_2=0$, $\rho_{12}=0.9$ and  $c_{11} = c_{22} = 0.6$.    The skewness of the simulated data is illustrated through the corresponding histograms.

\subsection{Results}

We  illustrate the saving of the pairwise  CL method in terms of computational burden.  All experiments were carried out on a 2.7 GHz  processor with 8 GB of memory and the estimation procedures were implemented coupling R functions and C routines. Table \ref{times_table} provides the  computational times (in seconds) in evaluating the weighted CL method, with cut-off distance equal to $d_{ij}=0.25, 0.5, 0.75, 1$ radians, for  all  $i,j=1,2$. These results show that CL  has a moderate computational cost even for large data sets. Indeed, the most demanding part in the  evaluation of the objective function  is the repeated numerical calculation of the bivariate Gaussian cdf.

\begin{table}
\caption{Time (in seconds) in evaluating CL method, with 0/1 weights, considering different number of observations and cut-off distances $d_{ij}=0.25,0.5,0.75,1$ (in radians).}
\label{times_table}
\begin{center}
\begin{tabular}{ccccccccc}
\hline\hline
\multicolumn{1}{c}{ }             &   \multicolumn{7}{c}{Number of observations }    \\  \cline{2-8}
               &  250         & 500      & 1000      &  2000   & 4000    &  8000  & 16000\\ \hline
$d_{ij}= 0.25$          &  0.003      & 0.007   &  0.016    &  0.068  & 0.254   &  0.954 &  3.753            \\  
$d_{ij}=0.5$           &  0.005      & 0.012   &  0.033    &  0.134  & 0.498   &  1.910 &  8.021            \\  
$d_{ij}=0.75$          &  0.007      & 0.017   &  0.050    &  0.205  & 0.796   &  3.020 &  12.024            \\  
$d_{ij}=1$          &  0.008      & 0.022   &  0.066    &  0.290  & 1.139   &  4.357 &  17.365            \\  \hline
\end{tabular}
\end{center}
\end{table}

We now study the statistical efficiency of the estimation method.  We consider  289  spatial sites in a grid on $\mathbb{S}^2$, which is generated  with 17 equispaced  longitude and  latitude points.  We use the latent correlation functions (\ref{exp_corr}) and (\ref{askey_corr}), with $\sigma_1^2=\sigma_2^2=1$, $\mu_1=\mu_2=0$, $c_{11}= 0.15$,  $c_{22}=0.25$ and $\rho_{12}=0.5$, under  the following choices for the skewness parameters:
\begin{itemize}
\item[{(A)}]   We set $\eta_1=1$ and $\eta_2=2$. In this case,  both components are right-skewed and the  collocated correlation coefficient between the components of the field is approximately $0.45$.
\item[{(B)}]  We set $\eta_1=1$ and $\eta_2=-2$. The first component of the field is right-skewed, whereas the second one is left-skewed,  and the collocated correlation coefficient between the components of the field  is approximately $0.19$. 
\end{itemize}
In total, we have four scenarios:
\begin{itemize}
\item \textbf{Scenario (I)}. Exponential model under choice (A).
\item \textbf{Scenario (II)}. Askey model under choice (A).
\item \textbf{Scenario (III)}. Exponential model under choice (B).
\item \textbf{Scenario (IV)}. Askey model under choice (B).
\end{itemize}

\begin{figure}
\centering{
\includegraphics[scale=0.42]{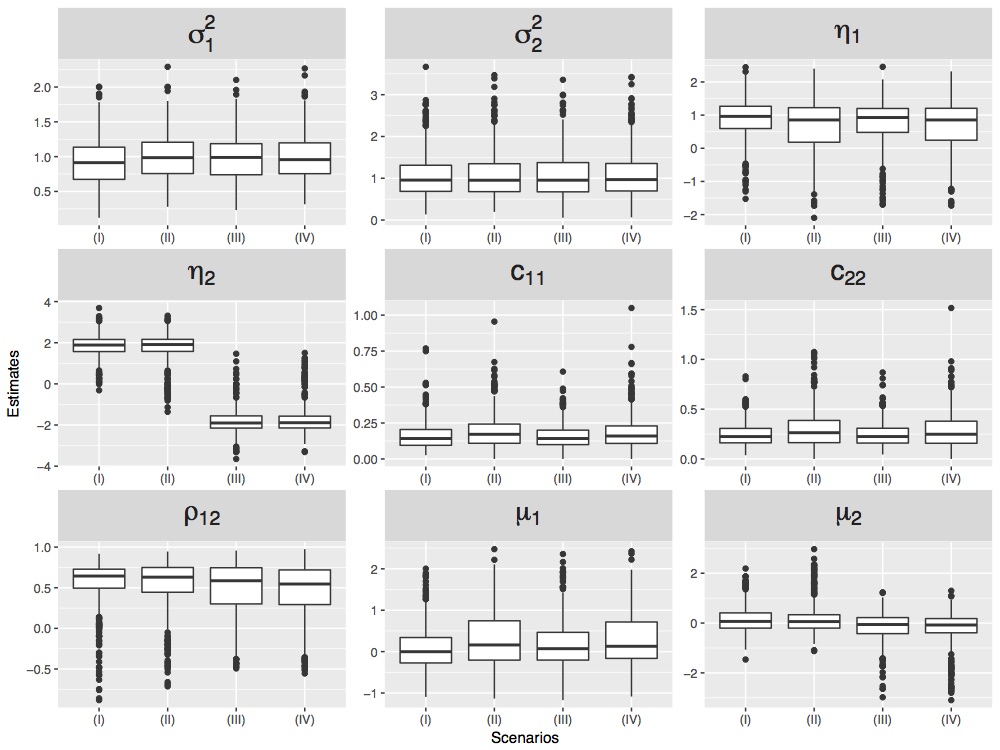} 
}
\caption{Boxplots of the CL estimates for the bivariate  skew-Gaussian RF,  under Scenarios (I)-(IV).  }
\label{boxplot}
\end{figure}

For each scenario, we simulate  500 independent  realizations from the bivariate skew-Gaussian RF. Then, we estimate the parameters using the weighted pairwise CL method. We set $d_{ij}=0.5$ radians for the 0/1 weights, for  all  $i,j=1,2$.   Figure \ref{boxplot}   reports the  boxplots of the CL estimates.     All studies show the effectiveness of our proposal. We have also  applied CL estimation to other parametric models, such as, the Cauchy and Wendland correlation functions. For each case considered the pairwise CL performs well.


\begin{figure}
\centering{
\includegraphics[scale=0.42]{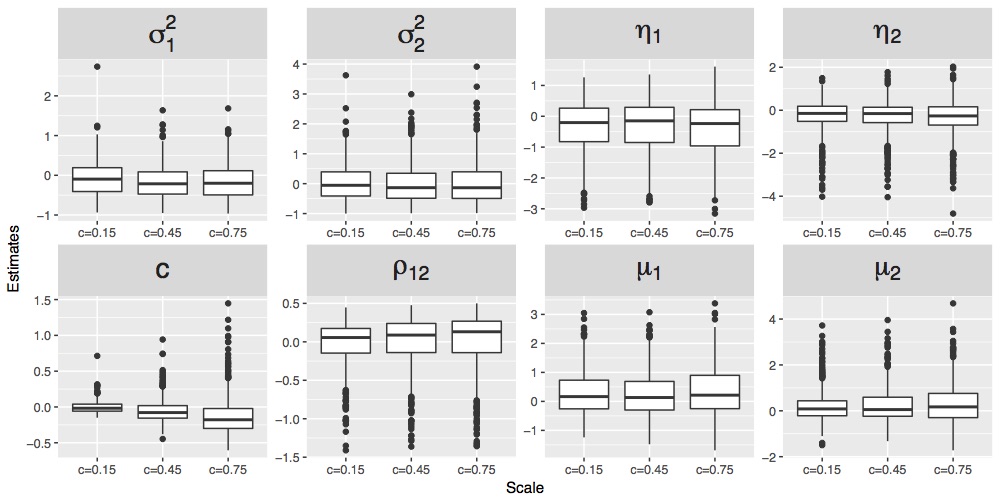}
} 
\caption{Centered boxplots of the CL estimates, for the bivariate  skew-Gaussian RF, using an Exponential latent correlation and different  scale parameters: $c=0.15, 0.45, 0.75$. }
\label{boxplot_scale}
\end{figure}

\begin{figure}
\centering{
\includegraphics[scale=0.42]{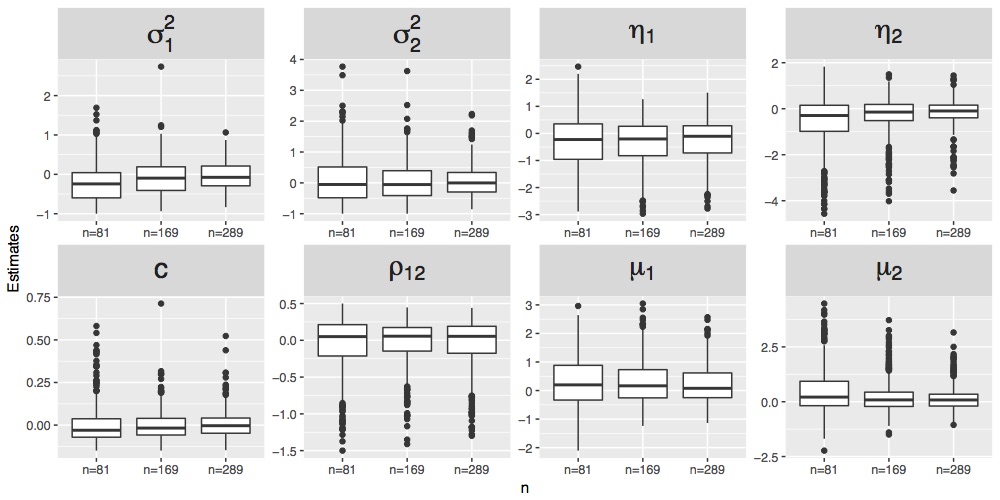}
} 
\caption{Centered boxplots of the CL estimates, for the bivariate  skew-Gaussian RF, using an Exponential latent correlation and different sample sizes: $n=81, 169, 289$.  }
\label{boxplot_n}
\end{figure}

Finally, we assess the performance of the pairwise CL method with increasing scale parameters $c_{ij}$ as well as increasing sample sizes $n$.  For simplicity, all the subsequent experiments consider a single scale parameter $c_{ij} = c$, for all $i,j=1,2$.  In this case, the parameter vector reduces to $\bm{\lambda} = (\sigma_1^2, \sigma_2^2, \eta_1, \eta_2, c, \rho_{12}, \mu_1, \mu_2)^\top$. We consider an Exponential latent correlation under the following parametric setting: $\sigma_1^2=\sigma_2^2=\eta_1=1$, $\eta_2=2$, $\mu_1=\mu_2=0$ and $\rho_{12}=0.5$.  Figure \ref{boxplot_scale} reports the centered boxplots of the  CL estimates in three different cases:  $c=0.15,0.45,0.75$. The increase of the parameter $c$ imply that the spatial dependence will be strengthened, and it produces biased estimates of $c$ and $\rho_{12}$. Our findings are consistent and add more evidence to the results reported in the previous literature  (\citealp{zhang2004inconsistent};  \citealp{Bevilacqua:Gaetan:Mateu:Porcu:2012}; \citealp{XU2016431}). On the other hand, we set $c=0.15$, and we consider increasing sample sizes: $n=81,169, 289$, in grids generated with 9, 13 and 17 equispaced longitude and latitude points, respectively.  As expected, more observations produce better estimations in terms of variability and bias.

\section{A bivariate data set}

We analyze a bivariate data set of Minimum (Variable 1) and Maximum (Variable 2)  surface air temperatures. The spatial variability of temperatures is crucial for modeling  hydrological and agricultural phenomena. These data outputs come  from the  Community Climate System Model (CCSM4.0) (see \citealp{doi:10.1175/2011JCLI4083.1})    provided  by  NCAR (National Center for Atmospheric Research) located at Boulder, CO, USA.  

\begin{figure}
\centering{
\includegraphics[height=7cm,width=7cm]{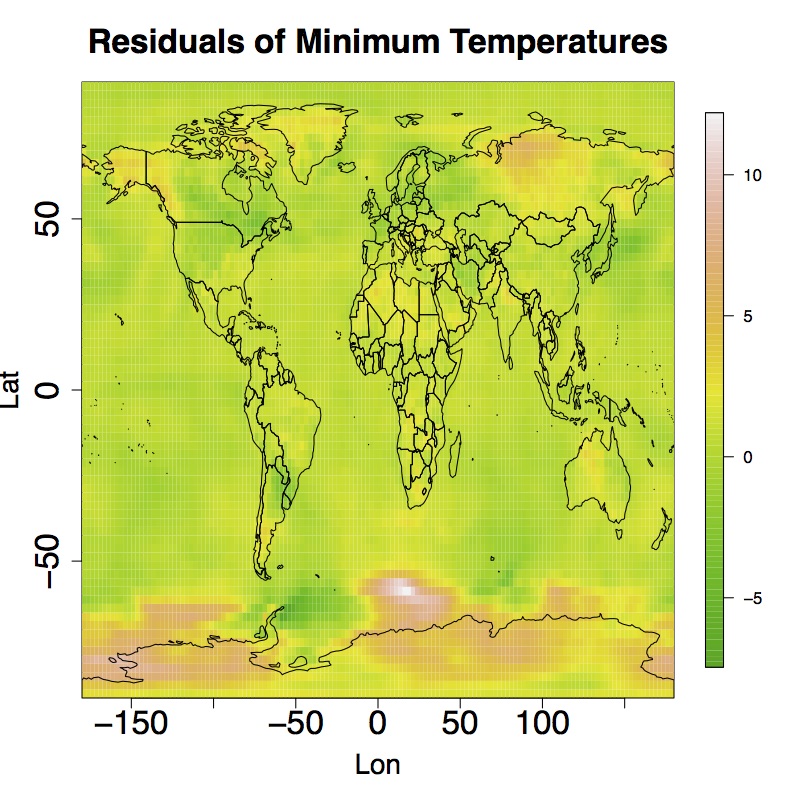} \includegraphics[height=7cm,width=7cm]{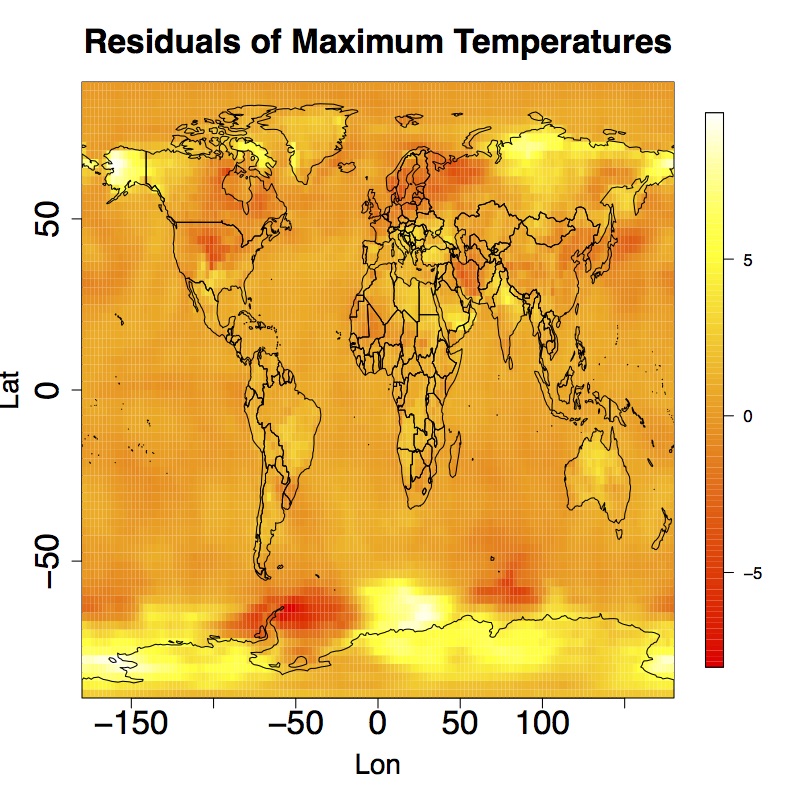}
} 
\caption{Residuals of the Minimum (left) and Maximum (right) surface air temperatures  in July of 2015.}
\label{real-data}
\end{figure}

We have monthly data  over a grid of  $2.5\times 2.5$  degrees of spatial resolution.  The unit for temperatures is Kelvin degrees.   We focus on  July of 2015 and we subtract the historical location-wise July average (considering the previous 50 years).  Figure \ref{real-data} depicts  the resulting residuals for the global data set. In order to ensure spherical isotropy, we only consider  locations with latitudes between $-30$ and $30$ degrees. The final data set consists of $3456$ observations per each variable.  These variables  are strongly correlated, since the empirical correlation  is 0.68.  The histogram of each variable reflects a certain degree of right skewness (see  Figure \ref{hist-data} below).  Thus, the residuals can be modeled approximately with our proposal, considering planet Earth as a sphere of radius 6378 kilometers.

We fit a bivariate  skew-Gaussian RF  with latent correlations of Exponential type. We use as benchmark a purely Gaussian model by taking   Equation (\ref{model_skew}) with $\eta_i=0$, for $i=1,2$.   We have considered the  parameterization introduced in the previous sections, so that, the skew-Gaussian model has 9 parameters, whereas the Gaussian model has 7 parameters. The CL estimation is carried out using only pairs of observations whose spatial distances are less than 1592.75 kilometers  (equivalent to $0.25$ radians on the unit sphere).  Table \ref{estimates_table} reports the CL  estimates for the skew-Gaussian and Gaussian models. The units of the scale parameters are kilometers.

\begin{figure}
\centering{
\includegraphics[height=5cm,width=5cm]{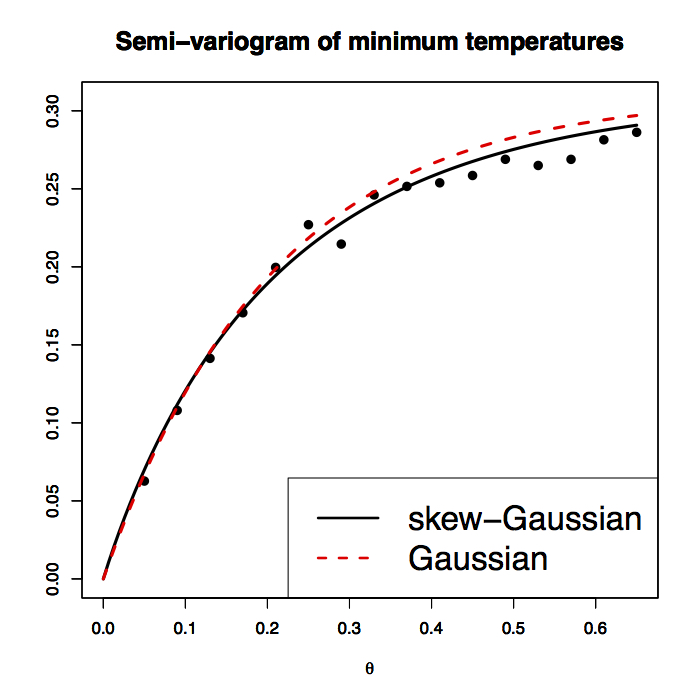} \includegraphics[height=5cm,width=5cm]{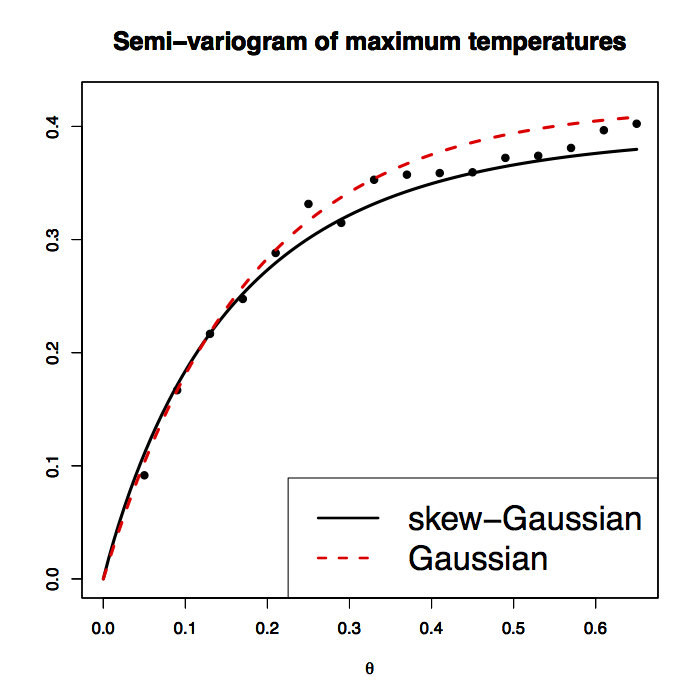}
\includegraphics[height=5cm,width=5cm]{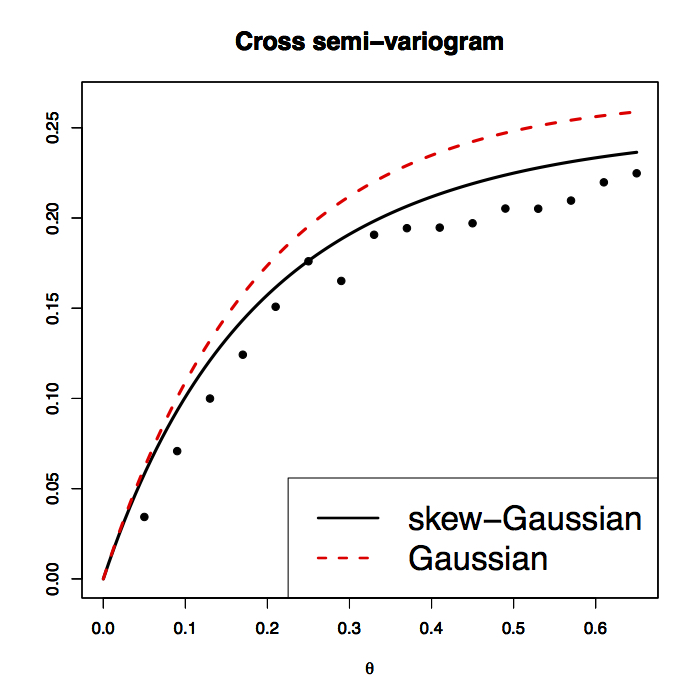}  
} 
\caption{Empirical semi-variograms versus fitted semi-variograms, using    Exponential latent correlations, for the skew-Gaussian (solid line) and Gaussian (dashed line) models.}
\label{variogram-data}
\end{figure}

The optimal values of the CL objective functions are given in Table \ref{pred_table}.   Note that the maximum  CL value under the  skew-Gaussian model is superior to the merely Gaussian model. It is clear that the incorporation of skewness produces improvements in goodness-of-fit.  Figure \ref{hist-data} shows the histograms of each variable and the fitted skew-Gaussian and Gaussian density functions. In Figure \ref{variogram-data}, the  marginal and cross empirical semi-variograms are compared to the theoretical models.

\begin{figure}
\centering{
\includegraphics[height=7cm,width=7cm]{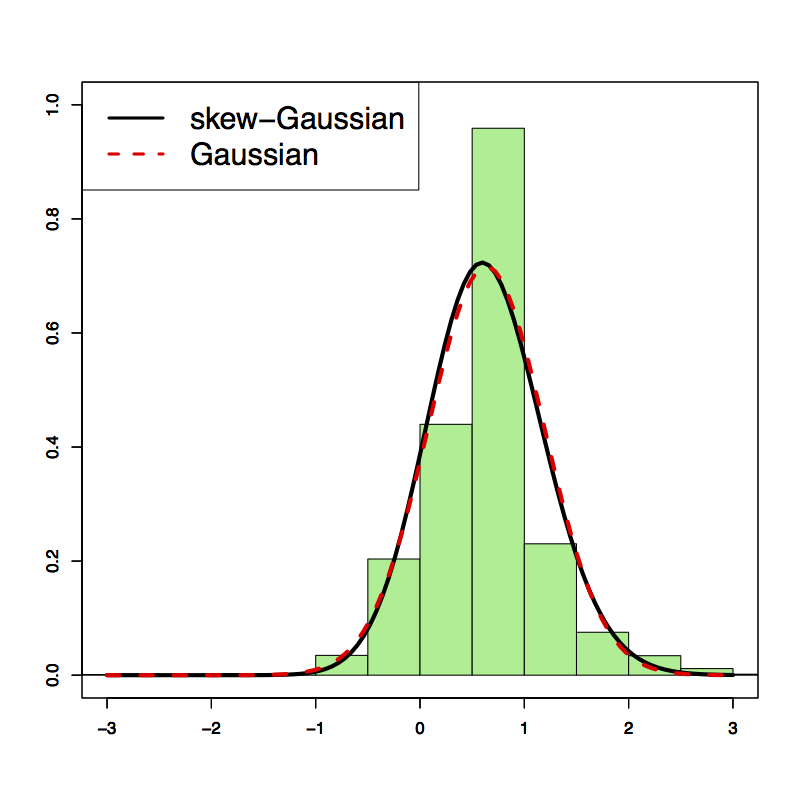} \includegraphics[height=7cm,width=7cm]{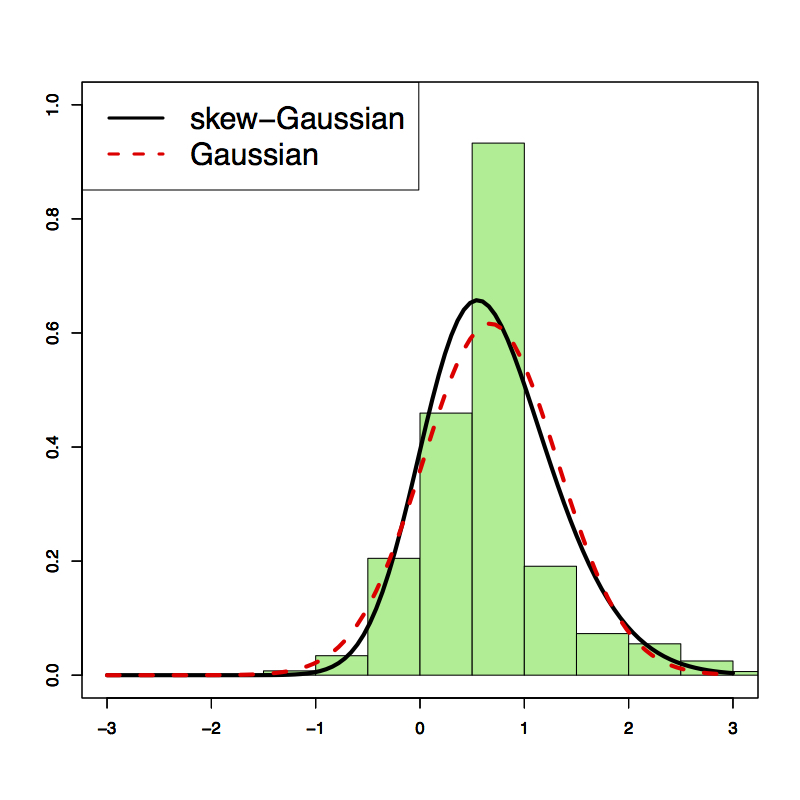}
} 
\caption{Histograms for the residuals of the Minimum (left) and Maximum (right) surface air temperatures, considering observations with latitudes between $-30$ and $30$ degrees, and the fitted skew-Gaussian (solid line) and Gaussian (dashed line) probability density functions.  }
\label{hist-data}
\end{figure}

  Finally, we compare   both models in terms of their predictive performance. Since the covariance structure of the skew-Gaussian field is known explicitly, we use  the classical best linear unbiased   predictor (cokriging), which is optimal for the Gaussian model,  in terms of mean squared error. However, it is not optimal for the skew-Gaussian RF.  In spite of this, we will show that the skew-Gaussian model provides better predictive results. We use a drop-one prediction strategy and quantify the discrepancy between the real and predicted values through the root mean squared prediction error (RMSPE)   
  $$ \text{RMSPE}= \sqrt{\frac{1}{2n}\sum_{i=1}^2 \sum_{k=1}^n (Z_i(\bm{s}_k) - \widehat{Z}_i(\bm{s}_k)  )^2} $$
   and the Log-score (LSCORE)
   $$  \text{LSCORE}=  \frac{1}{2n}\sum_{i=1}^2 \sum_{k=1}^n  \left[  \frac{\log(2\pi\widehat{\sigma}^2_i(\bm{s}_k))}{2}  + \frac{(Z_i(\bm{s}_k) - \widehat{Z}_i(\bm{s}_k)  )^2}{2 \widehat{\sigma}^2_i(\bm{s}_k)}   \right], $$
   where $n$ is the number of spatial locations, $\widehat{Z}_i(\bm{s}_k)$ is the drop-one prediction of $Z_i(\bm{s}_k)$ at location $\bm{s}_k$ and $\widehat{\sigma}^2_i(\bm{s}_k)$ is the drop-one prediction variance (see \citealp{zhang2010kriging}). Note  that the skew-Gaussian model generates better results since the mentioned indicators are smaller. In terms of RMSPE, the improvement in the prediction is approximately $3.1\%$.

\begin{table}
\caption{CL estimates for the skew-Gaussian and Gaussian RFs, using    Exponential latent correlations. Scale parameters are given in kilometers.}
\label{estimates_table}
\begin{center}
\begin{tabular}{cccccccccc}
\hline\hline
                     &  $\widehat{\sigma}_1^2$    &     $\widehat{\sigma_2}^2$      &    $\widehat{\eta}_1$    &    $\widehat{\eta}_2$     &     $\widehat{c}_{11}$
                     &     $\widehat{c}_{22}$       &        $\widehat{\rho}_{12}$       &       $\widehat{\mu}_1$     &      $\widehat{\mu}_2$  \\ \hline
skew-Gaussian  &    0.223   &  0.179  &  0.487 &   0.769   & 4995.9   &   4867.1   &  0.819  &   0.250  &  0.079     \\
Gaussian           &    0.310   &  0.419 &  -        &      -        &  3922.4 &  3393.1   &  0.743  &  0.635   &  0.674      	\\ \hline 
\end{tabular}
\end{center}
\end{table}

\begin{table}
\caption{Prediction performance of the skew-Gaussian and Gaussian RFs, using   Exponential latent correlations.}
\label{pred_table}
\begin{center}
\begin{tabular}{ccccc}
\hline\hline
                   &           $\#$ parameters &  CL &   \multicolumn{1}{c}{RMSPE}  &    \multicolumn{1}{c}{LSCORE}  \\  \hline
skew-Gaussian       &   9 &    $-1095912$       &   0.219             &     $-0.101$      \\   
Gaussian                &   7 &       $-1127334$     &      0.226    &    0.028  \\\hline
\end{tabular}
\end{center}
\end{table}

\section{Discussion}

Building models for non-Gaussian RFs has become a major challenge and more efforts should be devoted to such constructions. In particular, it seems that the main difficulties arise when trying to build models that are statistically identifiable. Another major problem, on the other hand, comes when building the finite dimensional distributions, which are analytically intractable in most cases. This paper has provided an approach that allows to avoid the identifiability problem in multivariate skew-Gaussian RFs, and that permits to implement a CL approach.  

We have shown that the pairwise CL method performs well  under different correlation structures and parametric settings. We believe that CL approaches can be adapted to other families of non-Gaussian fields, such as, Student's t or Laplace RFs, among others.   At the same time, the real data example illustrates that  the incorporation of skewness can produce significant improvements in terms of prediction in comparison to a Gaussian model.

Since the optimal predictor for the skew-Gaussian model, with respect to a squared error criterion,  is non-linear and   difficult to evaluate explicitly,  a  relevant research direction is the search for   methods that approximate this   predictor.      Indeed,  Monte Carlo methods are an appealing option \citep{zhang:el-shaarawi:2010}. However,  from a computational point of view, Monte Carlo samples are difficult to produce efficiently   and  such a method can be unfeasible  for large data sets.

\section*{Acknowledgments}
Alfredo Alegr\'ia is supported by \textit{Beca CONICYT-PCHA/Doctorado Nacional/2016-21160371}. Moreno Bevilacqua  is partially supported by \textit{Proyecto Fondecyt   1160280}.  Emilio Porcu is supported by \textit{Proyecto Fondecyt Regular number 1130647}.  Jorge Clarke is supported by \textit{Proyecto Fondecyt Post-Doctorado number 3150506.}

We also acknowledge the World Climate Research Programme's Working Group on Coupled Modelling, which is responsible for Coupled Model Intercomparison Project (CMIP).

 \appendix
 \section{Pairwise distributions for the multivariate skew-Gaussian RF}
 \label{apendice}

Before we state the proof of Proposition \ref{biv_log_lik},  we need the following property of quadratic forms. 

 \textbf{Lemma} Let $A$ and $B$ be two symmetric positive definite matrices of order $(n\times n)$ and $\bm{x},\bm{a}\in\mathbb{R}^n$. Then, we have the  following identity:
  \begin{equation} \label{identities}
  (\bm{a}-\bm{x})^\top A^{-1}(\bm{a}-\bm{x}) + \bm{x}^\top B^{-1}\bm{x} = (\bm{x}-\bm{c})^\top(A^{-1}+B^{-1})(\bm{x}-\bm{c})   + \bm{a}^\top(A+B)^{-1}\bm{a},
\end{equation}
 where $\bm{c}=(A^{-1}+B^{-1})^{-1}A^{-1}\bm{a}$.

\textbf{Proof of Proposition \ref{biv_log_lik}} Consider  $\bm{W}=(|X_{1}|,|X_{2}|)^\top$, $\bm{V}=(Y_{1},Y_{2})^\top$, $\bm{\mu}=(\mu_1,\mu_2)^\top$, $\bm{\eta}=(\eta_1,\eta_2)^\top$ and $\bm{\sigma} = (\sigma_1,\sigma_2)^\top$, where $(X_1,X_2)^\top \sim \mathcal{N}_2(\bm{0},\Omega(r^x))$ and $(Y_1,Y_2)^\top \sim \mathcal{N}_2(\bm{0},\Omega(r^y))$ are independent, with $\Omega(r)$ as defined in (\ref{corr_matrix}).
 Let $\bm{Z}=(Z_{1},Z_{2})^\top$ defined through $$\bm{Z}=\bm{\mu}+\bm{\eta} \circ \bm{W}+ \bm{\sigma} \circ \bm{V},$$ where $\circ$ denotes the Hadamard product.
Therefore, the joint probability density function of $\bm{Z}$ is given by
\begin{equation}
 f_{\bm{Z}}(\bm{z})= \displaystyle  \int_{\mathbb{R}^{2}_{+}} f_{\bm{Z}|(\bm{W}=\bm{w})}(\bm{z}|\bm{w})f_{\bm{W}}(\bm{w}) \text{d}\bm{w},
  \label{biv_skew}
\end{equation}

Here, $f_{\bm{W}}(\bm{w})$ is the pdf of the random vector $\bm{W}$ and  $\bm{w} = (w_1,w_2)^\top$. Note that the cdf   of the random vector $\bm{W}$  can be written as
\begin{equation*}
F_{\bm{W}}(w_{1},w_{2})  =  \Phi_{2}(w_{1},w_{2};\Omega(r^x)) - \Phi_{2}(-w_{1},w_{2};\Omega(r^x)) 
- \Phi_{2}(w_{1},-w_{2};\Omega(r^x)) +\Phi_{2}(-w_{1},-w_{2};\Omega(r^x))
\end{equation*}
Then, we can  obtain  the pdf of $\bm{W}$,
\begin{eqnarray*}
f_{\bm{W}}(w_{1},w_{2}) & = &  2  \bigg( \phi_{2}(w_{1},w_{2};\Omega(r^x)) + \phi_{2}(-w_{1},w_{2};\Omega(r^x)) \bigg) ,\\
                  & = & 2  \bigg( \phi_{2}(w_{1},w_{2};\Omega(r^x)) + \phi_{2}(w_{1},w_{2};\Omega(-r^x)) \bigg).
\end{eqnarray*}
On the other hand, $f_{\bm{Z}|\bm{W}=\bm{w}}(\bm{z}|\bm{w})$ is the pdf of the random vector $\bm{Z}|(\bm{W}=\bm{w})\sim \mathcal{N}_{2}(\bm{\mu}+\bm{\eta}\circ \bm{w};\Omega_2)$, with $$\Omega_2 = \begin{pmatrix}  \sigma_1^2  &  \sigma_1\sigma_2  \\ \sigma_1\sigma_2 &  \sigma_2^2\end{pmatrix} \circ \Omega(r^y).$$
Therefore, evaluation of the integral (\ref{biv_skew}) requires the characterization of an integral of the form 
\begin{eqnarray*}
 \mathcal{I} &  =  & \displaystyle  \int_{\mathbb{R}^{2}_{+}} \phi_{2}(\bm{z}-\bm{\mu}- \bm{\eta}\circ \bm{w}; \Omega_{2}) \phi_{2}(\bm{w}; \Omega_{1}) \text{d}\bm{w},\\
\end{eqnarray*}
with $\Omega_{1}$ being $\Omega(r^x)$ or $\Omega(-r^x)$. Moreover, $\mathcal{I}$ can be written as
\begin{equation*}
 \mathcal{I} =  \displaystyle  |\Upsilon|  \int_{\mathbb{R}^{2}_{+}}  \phi_{2}(\Upsilon(\bm{z}-\bm{\mu}) -  \bm{w};\Upsilon \Omega_{2}\Upsilon) \phi_{2}(\bm{w}; \Omega_{1}) \text{d}\bm{w}
\end{equation*}
where $\Upsilon = \text{diag}\{1/\eta_1, 1/\eta_2\}$.
Thus, using  Equation (\ref{identities}) we have
\begin{eqnarray*}
\mathcal{I} &  = &  |\Upsilon| \displaystyle  \phi_{2}\bigg(\Upsilon(\bm{z}-\bm{\mu}); \Upsilon \Omega_2 \Upsilon +\Omega_{1}\bigg)  \int_{\mathbb{R}^{2}_{+}}   \phi_{2}\bigg( \bm{w}-L;           ( [\Upsilon \Omega_2 \Upsilon]^{-1} + \Omega_1^{-1} )^{-1} \bigg) \nonumber \text{d}\bm{w},\\
 &  =    &  |\Upsilon|   \phi_{2}\bigg(\Upsilon(\bm{z}-\bm{\mu}); \Upsilon \Omega_2 \Upsilon +\Omega_{1}\bigg)
 \Phi_{2}  ( L ;  ( [\Upsilon \Omega_2 \Upsilon]^{-1} + \Omega_1^{-1} )^{-1}    ) \nonumber\\
  &  =    &   \phi_{2}\bigg(\bm{z}-\bm{\mu}; \Omega_{2}+\Upsilon^{-1} \Omega_1 \Upsilon^{-1}\bigg)
 \Phi_{2}  (L;  ( [\Upsilon \Omega_2 \Upsilon]^{-1} + \Omega_1^{-1} )^{-1}    ) \nonumber
\end{eqnarray*}
 where $L = \bigg[  I_2 + \Upsilon \Omega_2 \Upsilon \Omega_1^{-1}   \bigg]^{-1} \Upsilon (\bm{z}-\bm{\mu})$.

 \section{Multivariate Askey model on the sphere}
 \label{apendiceB}

We consider a multivariate Askey model for the sphere $\mathbb{S}^2$, defined according to 
$$ \rho_{ij } r (\theta; c_{ij}) = \rho_{ij} \left ( 1- \frac{\theta}{c_{ij}} \right )_{+}^{4}, \qquad \theta \in [0,\pi], \quad i,j=1,\hdots,m, $$
which has been used through Sections 4 and 5, considering $c_{ij}=(c_{ii}+c_{jj})/2$. We claim that such a model is positive definite  using  the following mixture  (see \citealp{Daley:2014})
$$ r (\theta; c_{ij}) \propto \int_{0}^{\infty} r(\theta; \xi) \xi^2 r(\xi;c_{ij}) {\rm d} \xi,   \qquad \theta \in [0,\pi].$$
The results of \cite{gneiting2013} for the univariate  Askey model, coupled with the conditions developed by \cite{Daley:2014}  complete our assertion.



\section*{\refname}
 \bibliographystyle{elsarticle-harv} 
\bibliography{mybib}

\end{document}